\documentclass[12pt]{article}
\usepackage{graphicx}
\usepackage{geometry}
\usepackage{relsize}
\author{\relsize{+1}Vibeke Libby\\\relsize{-1} Woodside, CA 94062\\ \\ \relsize{-1} \sl www.vibekelibby.com\\  \relsize{-1}\sl agilemath@gmail.com}
\title{\bf A Fast Algorithm for Determining the Existence and Value of Integer Roots of $N$}
\begin{document} 
\maketitle
\begin{abstract}
We show that all perfect odd integer squares not divisible by $3$, can be usefully written as $\sqrt{N}=a+18p$, where the constant $a$ is determined by the basic properties of $N$. The equation can be solved deterministically by an efficient four step algorithm that is solely based on integer arithmetic. There is no required multiplication or division by multiple digit integers, nor does the algorithm need a seed value. It finds the integer $p$ when ${N}$ is a perfect square, and certifies $N$ as a non-square when the algorithm terminates without a solution. The number of iterations scales approximately as $log_9(\sqrt{N}/2)$ for square roots. The paper also outlines how one of the methods discussed for squares can be extended to finding an arbitrary root of $N$. Finally, we present a rule that distinguishes products of twin primes from squares.
\end{abstract}

\section{Problem and Summary}
In cryptography and broader applications, such as the general factorization of integers, it is often of interest to determine quickly if an integer $N$ is a true square or is a product of two closely spaced integers, such as twin primes. Numerous non-deterministic methods of varying efficiencies are currently applied to the problem of finding the square root of an integer. These methods include digit-by-digit calculations, Taylor series, Newton's method {\it ref [1]}, Pell's equation {\it ref [2]}, and even sophisticated hardware implementations {\it ref [3]}. For very large integers, these methods, may not yield rapid conclusive results. In contrast, the algorithm presented here provides a conclusive answer after a specific number of steps, determined by the size of the $N$. Throughout this presentation, all variables and computations are based on positive integer arithmetic alone.

The number of required iterations to find the square root of $N$ is approximately $log_9(\sqrt{N}/2) - 3/2$. The same number of iterations is required to certify that $N$ is a non-square.

\section{Properties of N}
Suppose the object is to find the square root of a general integer $M$. The presented algorithm works by focusing on $N$, the odd factor of $M$ that is not divisible by 2 or 3 to any power.

\begin{equation}
M = 2^k*3^l*N
\end{equation}
Clearly, if either $k$ or $l$ is odd, $M$ is not a perfect square. If $k$ and $l$ are both even, then $M$ is square if $N$ is square. Since factors of $2$ and $3$ can quickly be identified in even large numbers, we can choose a representation for $N$ which do not contain these divisors.

One such representation is the reduced residue system $U_{18}$ with the elements \\$\{1, 5, 7, 11, 13, 17\}$ or $\{\pm1, \pm5, \pm7\}$ modulo $18$, where each element $[r]$ is a residue class in $U_{18}$. By definition, $ref [4]$, all $[r]$ are relatively prime to $18$ from which it follows that g.c.d. $(r,18) = 1$ and the congruence $r*N\equiv 1\pmod{18}$ has a solution. It also happens to be true that all $[r]$ in $U_{18}$ are relatively prime to each other, which is not generally true for reduced residue systems; but is true, for example, for systems modulo $3$, $6$, $8$, $10$, and $18$. 

Since the integer $18$ only contains the divisors $2$ and $3$, $U_{18}$ contains no even numbers and no odd numbers divisible by $3$. Therefore, $U_{18}$ is an ideal representation of all possible $N$-values because every $N$ defined in equation 1 can be $\it uniquely$ assigned to one of the six elements \{1, 5, 7, 11, 13, 17\}. According to Euler's theorem, for $\forall r$ relatively prime to $18$, $r^{6n} \equiv 1 (mod18)$. See figure 1.

To test for divisibility by 3 one traditionally computes the sum-of-digits (SOD) of $N$ repeatedly until a number between one and nine results, often called the Digit Root (DR). If DR is 3, 6, or 9 then $N$ is divisible by three. This test is equivalent to computing $N(mod 9)$ where the latter calculation involves division by 9 as opposed to summation of digits. 

Here, we extend this test to cover $N(mod18)$ replacing division by $18$ with the addition of the digits in $N$. We will show that a mapping exists between the elements in $U_{18}$, $N(mod18)$, and $DR(N)$, $N(mod9)$, of the form:

\begin{equation}
[N(mod18)](mod9) = N(mod9)
\end{equation}

$Proof:$\\
For $c< 18$ write N(mod18) as $N = c + 18X = c + 9(2X)$ and for $d < 9$ write N(mod9) as $N = d + 9Y$. It follows that $(c + 9(2X))(mod9) = c(mod9)$ and $(d + 9Y)(mod9) = d(mod9) = d$. Therefore, the equality in equation 2 is satisfied when:
\begin{equation}
c(mod9) = d
\end{equation}

This condition is uniquely met for all $N$ in $U_{18}$; and therefore the DR method can be applied to the elements in $U_{18}$ as long as it is recognized that we have to use c  = 11 for d = 2, c = 13 for d = 4, and c = 17 for d = 8. For 1, 5, 7: c = d.
To find which element an arbitrary $N$ belongs to, we apply the condition in equation 3. By doing so, we have substituted division by $18$ with addition. For example, if $N = 512346251$, SOD(N) = 29 with DR(29) = 2, then N belongs to element $11$, written as $[11]$ or 11 + 18X.

\section{Determining solutions for $a$ for known $N$}
Writing any integer in $U_{18}$ as $a + 18p$, with $a$ is uniquely determined by DR(N), any square in $U_{18}$ can be written as:
\begin{equation}
N  = (a + 18p)^2 = a^2 + 36ap + 18^2p^2
\end{equation}

From this expression $DR(N)$ equals $DR(a^2)$ because the remaining terms contain the factor $18$, which does not affect the $DR$. Since $a$ is limited to the values 1, 5, 7, 11, 13, 17, the squares in $U_{18}$ are therefore found only in: $[1]$, $[7]$, $[13]$. This simplification comes with a computational cost for squares because two $a$-values are possible for each element: 
\begin{itemize}
\item $[1]$ for a = 1, 17 
\item $[7]$ for a = 5, 13 
\item $[13]$ for a = 7, 11 
\end{itemize}
Therefore, two parallel tests are needed for each square root calculation.
\begin{figure}[ht]\label{fig:in1}
\begin{center}
\begin{tabular}{|c|c|c|c|c|c|c|}
\hline
N(mod 18)=[ ] & {} & {} & {} & {} & {} & {}\\ 
\hline
\hline
$a$ & $[a]^{6n}$&$[a]^{1+6n}$& $[a]^{2+6n}$ & $[a]^{3+6n}$ & $[a]^{4+6n}$ & $[a]^{5+6n}$ \\
\hline
$1$ & [1] & [1] & [1] & [1] & [1] & [1]\\
\hline
$5$ & [1] & [5] & [7] & [17] & [13] & [11] \\
\hline
$7$ & [1] & [7] & [13] & [1] & [7] & [13]\\
\hline
$11$ & [1] & [11] & [13] & [17] & [7] & [5]\\
\hline
$13$ & [1] & [13] & [7] & [1] & [13] & [7]\\
\hline
$17$ & [1] & [17] & [1] & [17] & [1] & [17]\\
\hline
\end{tabular}
\caption{\it The complete $mod 18$ power multiplication table for $U_{18}$. Example: $a = 5$ to the $4^{th}$ power belongs to residue class [13] because $5^4 = 625$ and $N$ is therefore congruent to $13 (mod 18)$. Conversely, $N$ can't be a power of $2 + 6n$, and specifically a square, if $N mod 18$ has a value other than 1, 7, or 13.}
\end{center}
\end{figure}

As argued above for squares, the value of $N(mod18)$ is required to determine the possible values of $a$ for any root. From the first, third, and fifth column in figure 1, all even powers belong to the residue classes $[1]$, $[7]$, $[13]$. Also note that six tests for the constant $a$ are required for powers of $6n$, three tests for powers of form $3+6n$, two tests for powers of $2+6n$ and $4+6n$, and only one test each for powers of $1+6n$ and $5+6n$. All such tests are independent and can be performed in parallel. 

\section{Algorithm Description - Squares}
The algorithm is divided into three parts:
\begin{enumerate}
\item Setup
\item Initialization
\item Loop
\end{enumerate}
\subsection{Setup}
With $N$ and $a$ given rewrite condition 4 as:
\begin{equation}
\frac{N_0 - ap}{9} = p^2
\end{equation}
Here $N_0 = \frac{N - a^2}{36}$ is guaranteed integer when $(DR(N),a)$ assumes any of the six combinations $(1,1), (1,17), (7,5), (7,13), (13,7), (13,11)$, see figure 1. Therefore $N$ is square with a root of the form $(a + 18p)$ if there exists an integer solution for $p$ satisfying condition 5. We will iteratively test for the existence of such $p$ and for this purpose define the $\it unique$ representation of $p$ as:
\begin{equation}
p = \sum_{i=0}^{i-max} 9^i*b_{i+1}
\end{equation}
with $1 \leq b_{i+1} \leq 9$ and where i-max is the smallest i-value for which the termination criterion:
\begin{equation}
\frac{N_0 - ap_i}{9} \leq p_i^2
\end{equation}
is met. When the equality sign is valid, then $N$ is square; and when the reverse is true a solution may still exist. In this case the iterations must continue until the termination criterion is met.

\begin{equation}
p_i = b_1 + 9b_2 + 9^2b_3 + ... + 9^ib_{i+1} = p_{i-1} + 9^ib_{i+1}
\end{equation}

\subsection{Initialization}
Set $i = 0$, $b_0 = 0$, and $p_{-1} = 0$.

\begin{figure}[ht]\label{fig:in2}
\begin{center}
\begin{tabular}{|c|c|c|c|c|c|}
\hline
\bf{Initialization} & $\bf{N_0}$ & $\bf{frac_0}$ & $\bf{b_0}$ & $\bf{p_0}$ & $\bf{f_0}$ \\
\hline
{i = 0} & \relsize{+1.5} $\frac{N-a^2}{36}$ & \relsize{+1.5} $\frac{N_0 - a*b_1}{9}$ & 0 & $b_1$ & $b_1*p_0$  \\
\hline
\end{tabular}
\caption{\it $Init$ shows the algorithm initial conditions, where $b_1$ is determined from the condition that $frac_0$ be integer.}
\end{center}
\end{figure}
It is useful to define the quantity $f_0 = b_1 * b_1 = p_0 * b_1$ and call the left hand side of equation 7 $frac_0$.

\subsection{Loop}
\subsubsection{i = 0}
Using equation 7 and the initial values from figure 2 condition 5 becomes an evaluation of:
\begin{equation}
\frac{N_0 - ab_1}{9} == b_1 * b_1
\end{equation}
We know from the unique representation of $p$ that there exists one and only one $b_1$ value for which $frac_0$ is integer and that $b_1$ must satisfy $DR(N_0) = DR(a*b_1)$. This condition can be implemented in a nine-entry look-up table. With $b_1$ now known, condition 7 can be evaluated. The evaluation terminates when $frac_0 < f_0$, continues for $frac_0 > f_0$, and declares a solution for $frac_0 = f_0$.

\subsubsection{i = 1}
In the next step $p_1 = b_1 + 9b_2$. Inserting $p_1$ and $p_1^2$ into condition 5 yields:
\begin{equation}
\frac{frac_0 - f_0 - ab_2}{9} = b_2(2b_1 + 9b_2)
\end{equation}
Setting $frac_0 - f_0 = N_1$ and realizing $2b_1 + 9b_2 = p_0 + b_1 + 9b_2$ yields a condition similar to 9; but now with $b_2$ as the unknown quantity: 
\begin{equation}
\frac{N_1 - ab_2}{9} == b_2(p_0 + p_1)
\end{equation} 
As before, $b_2$ is found from the integer condition of the left hand side, $frac_1$. As $i$ increases, the number of terms on the right hand side, $f_i$, increases as well. However, since each $f_i$ can be expressed in terms of the previous and current value of $p_i$, there exists a general recursion relation for $f_i$ as shown in figure 3. 

\subsubsection{General Case: i $\geq 1$}
One can express the i'th value, for $i > 0$, of $N_i$ as:
\begin{equation}
N_i = {\frac{N_{i-1} - a*b_i - 9^i*b_i^2 - 2*b_i*\sum_{k=0}^{k=i-1}9^k*b_k} {9}}
\end{equation}
and
\begin{equation}
\frac{N_i-a*b_{i+1}}{9} =  9^i*b_{i+1}^2 + 2*b_{i+1}*\sum_{k=1}^{k=i}9^{k-1}*b_k 
\end{equation}
Using the definition of $p_i$, equation 8, the right hand side of equation 13 becomes:
\begin{equation}
f_i = b_{i+1}*(9^i*b_{i+1} + 2*p_{i-1}) = b_{i+1}*(p_i + p_{i-1})
\end{equation}
The algorithm logic flow diagram is shown in figure 4.

\begin{figure}[ht]\label{fig:in3}
\begin{center}
\begin{tabular}{|c|c|c|c|c|c|}
\hline
\bf{i} & $\bf{N_i}$ & $\bf{frac_i}$ & $\bf{b_i}$ & $\bf{p_i}$ & $\bf{f_i}$  \\
\hline
$i>0$ &  $frac_{i-1}-f_{i-1}$ & \relsize{+1.5} $\frac{N_i-a*b_{i+1}}{9}$ &  $b_{i+1}$ &  $p_{i-1}+9^i*b_{i+1}$ &  $b_{i+1}*(p_{i-1}+p_i)$  \\
\hline
\end{tabular}
\caption{\it Illustration of the general recursion relation for all the variables as a function of $i$. Reading the tables from left to right gives the order of computation. The quantity $b_{i+1}$ can assume the values 1 through 9 and is uniquely determined from the integer condition: $\relsize{+1.5} {\frac{N_i -a*b_{i+1}}{9}}$.}
\end{center}
\end{figure}

\begin{figure}
\begin{center}
	\includegraphics{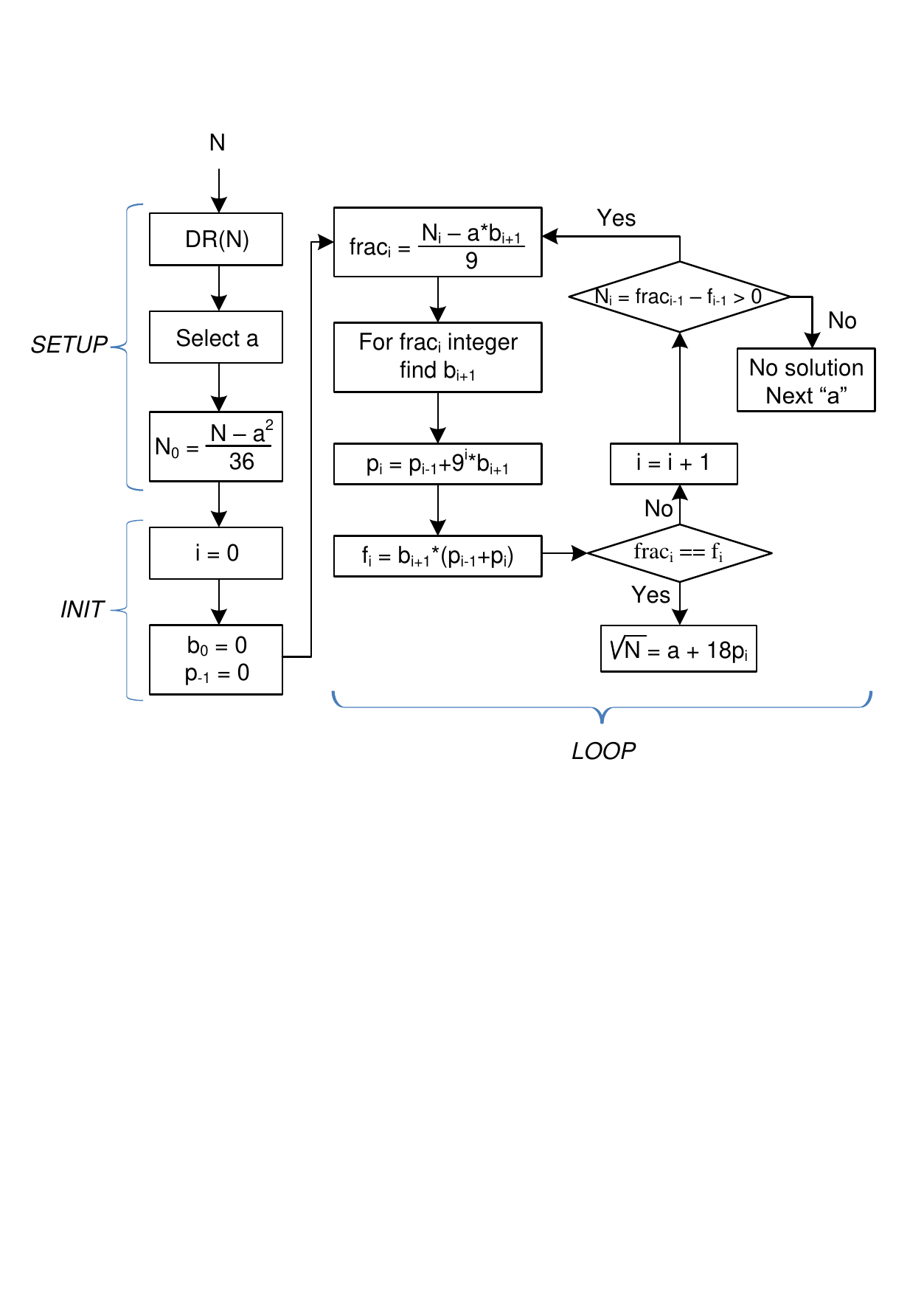}
	\caption{\it Logic flow diagram for finding the square root of an odd integer N not divisible by 3.}
\label{fig4:LibbySquare3}
\end{center}
\end{figure}

\section{Positive and Negative Certification}
Here we show that $N$ is a perfect square if and only if the condition $frac_i = f_i$ is satisfied after $i$ steps. Conversely, $N$ is not a square if $frac_j < f_j$ without having satisfied the condition for squareness above for $i < j$. Of general concern could be that the algorithm would map different $N$ values, squares and non-squares, to the same $p$. As we shall see, this is not the case.

Combining the square condition: $N = (a + 18p)^2$ with equation 6 yields the test condition for squareness:
\begin{equation}
N == (a + 18(b_1 + 9b_2 + ... + 9^ib_{i+1}))^2
\end{equation}
As discussed in section 4, equation 15 is equivalent to evaluating:
\begin{equation}
N_0 - a(b_1 + 9b_2 + ... + 9^ib_{i+1}) == 9(b_1 + 9b_2 + ... + 9^ib_{i+1})^2
\end{equation}
Section 4 also showed how each $1\leq b_i \leq 9$ can be uniquely determined for each step in $i$, with $frac_i$ decreasing and $f_i$ increasing with $i$ until condition $frac_i \leq f_i$ is satisfied or superseded. At step $i$ all values of $b$ are uniquely determined in the expansion of $p$; and because the representation of $p$ is a unique expansion of $b_i$ mod 9, we see the validity of the perfect square condition stated above. If $N$ tests to be square, its root is also unique per equation 15.

\section{Generalization to Cube and Higher Roots}
For cubes, $N_0$ is determined from $N = (a + 18p)^3$ as:
\begin{equation}
N_0 =\frac{N-a^3}{3*18} = a^2p + 18ap^2 + 6*18p^3
\end{equation}

Recall, from figure 1 that three parallel tests for $a$ have to be conducted for this case.
The first equation to be evaluated becomes:

\begin{equation}
\frac{N_0 - a^2b_1}{9} = 2ab_1^2 + 12b_1^3
\end{equation}

Although two divides both sides, we will maintain the same nomenclature as for squares. In this case $frac_0 = (N_0-a^2b_1)/9$ and $f_0 = 2b_1^2(a + 6b_1)$. Proceeding as for squares $N_1 = frac_0 - f_0$ becomes:
\begin{equation}
N_1 = \frac{N_0 - a^2*b_1 - 9*f_0}{9}
\end{equation}
We return to equation 18 and insert the next value for $p$; $p_1 = b_1 + 9b_2$ yielding:
\begin{equation}
\frac{N_0 - a^2(b_1 +9b_2)}{9} == 2a(b_1 + 9b_2)^2 + 12(b_1 + 9b_2)^3
\end{equation}
As in the quadratic case the termination criterion is met when: $frac_i \leq f_i$.

Observe that the same equations and steps applies to cubes as for squares. The difference between the two methods lies in the existence of a  recursion relation for squares. Perhaps recursion relations can be found for higher powers. If not, the above step-by-step approach can be applied to all powers. Although the formulae for a given power contain increasingly more terms as $i$ increases, they only have to be derived once per root.

\section {Number of Iterations Required}
The smallest possible value of $N$ of the form $(a + 18p)*(a + 18p)$ is $N = 361$ for $a = 1$ and $p = 1$. Smaller squares are possible for $p = 0$ but are considered trivial in this context since the focus here lies in computing roots of very large integers. In order to find an approximate formula for the number of iterations, $i$, we recall that $frac_i$ decreases by a factor of $9$ per step, whereas $f_i$ increases by a factor of $9$ per step, and that the algorithm terminates when the two quantities are equal or $frac_i < f_i$. Therefore, the relative size of $frac_i$ and $f_i$ after $i$ steps are:
\begin{equation}
\frac{frac_0}{9^i} = f_0 * 9^i
\end{equation}

From equation 5 and the conservative approximation of $a = 0$, we find $frac_0 = N/(9^2*2^2)$. Setting $f_0 = 9$ conservatively in equation 20 yields the desired expression between $N$ and the number of iterations $i$, $\frac{N}{4} = 9^{2i+3}$. This can be rewritten as:

\begin{equation}
i = log_9(\frac{\sqrt{N}}{2})-\frac{3}{2}
\end{equation}
Expression 21 represents the number of iterations required for each $a$ value and therefore doubles when two tests have to be performed. Recall from equation 1 that, $N$ does not contain any factors of two and three. If $N$ is not a square, the same number of iterations are required.

\section{Numerical Examples}
\subsection{Square Root Numerical Example}
As a point example consider $\sqrt{N}$ for $N = 1429822969 = 37813^2$ with $DR(N) = 7$. The allowed values for $a$ in this case are $13$ and $5$. We will test for $a = 13$ first and see how the algorithm works when a solution is found, in figure 5, and then show what happens when there is no solution by testing for $a=5$, figure 6. 

\begin{figure}[ht]\label{fig:in5}
\begin{center}
\begin{tabular}{|l|l|c|c|c|c|c|c|}
\hline
{i} & {a = 13} & $DR(N_i)$ & ${b_i}$ & $b_{i+1}$ & $p_i$ & $frac_i$ & $f_i$ \\
\hline
\hline
{} & $N = 1429822969$ & {7} &{} & {} & {} & {} & {} \\
\hline
{0} & $N_0 = 39717300$ & {3} & {0} & 3 & 3 & {4413029} & $9$ \\
\hline
{1} & $N_1=4413020$ & {5} & {3} & {8} & 75 & $490324$ & $624$ \\
\hline
{2} & $N_2=489700$ & {1} & {8} &{7} & 642 & $54401$ & $5019$ \\
\hline
{3} & $N_3=49382$ & {8} & {7} & {2} & 2100 & $5484$ & $5484$ \\
\hline
{4} & $N_4 = 0$ & {} & {2} & {} & {} & {} & {} \\
\hline
\end{tabular}
\caption{\it Illustration on how the square root of $N = 1429822969$ can be found in four iterations. The algorithm terminates with a solution when $frac_i = f_i$. Thus $p = 2100 = (3+9*8+9^2*7+9^3*2)$ with solution $\sqrt{N} = 13+18*2100 = 37813$.}
\end{center}
\end{figure}

\begin{figure}[ht]\label{fig:in6}
\begin{center}
\begin{tabular}{|l|l|c|c|c|c|c|c|}
\hline
{i} & {a=5} & $DR(N_i)$ & $b_i$ & $b_{i+1}$ & $p_i$ & $frac_i$ & $f_i$  \\
\hline
\hline
{} & $N = 1429822969$ & {} & {} & {} & {} & {} & {} \\
\hline
{0} & $N_0 = 39717304$ & 7& {0} & {5} & 5 & $4413031$ & $25$ \\
\hline
{1} & $N_1=4413006$ & 9 & {5} & {9} & 86 & $490329$ & $819$ \\
\hline
{2} & $N_2=489510$ & 9 & {9} & {9} & 815 & $54385$ & $8109$ \\
\hline
{3} & $N_3=46276$ & 7 & {9} & {5} & 4460 & $5139$ & $26375$ \\
\hline
{4} & $N_4 = -21236$& {} & {5}& {} & {} & {} & {}\\
\hline
\end{tabular}
\caption{\it Illustration on how a square root test of a ten-digit number $N$ is concluded in four steps when $N$ is not square. The algorithm terminates when $frac_i < f_i$ which occurs for $i = 3$.}
\end{center}
\end{figure}

\subsection {Cube Root Numerical Example}
The following example shows the algorithm steps for finding the cube root of a thirteen digit number $N = 6,177,847,762,549$. Since $DR(N) = 1$, $a$ can take on the values 1, 7, and 13. We will test $N$ for $a=7$ since this is the solution in this example.
The cubic root in this example is not computed by using a recursion relation but by step-by-step calculation of $frac_i$ and $f_i$.
\begin{figure}[h]\label{fig:in7}
\begin{center}
\begin{tabular}{|l|l|c|c|c|c|}
\hline
{i} & {a = 7} & $DR(N_i)$ & $b_{i+1}$  & $frac_i$ & $f_i$ \\
\hline
\hline
{} & $N = 6177847762549$ & {1}  & {} & {} & {}\\
\hline
{0} & $N_0 = 114404588189$ & 8 & {2} & $12711620899$ & $152$ \\
\hline
{1} & $N_1=12711620747$ & 2 & {5}  & $1412402278$ & $141850$ \\
\hline
{2} & $N_2=1412260428$ & 3 &{3} & $156917809$ & $3611958$ \\
\hline
{3} & $N_3=153305851$ & 4 & {1} &  $17033978$ & $17033978$ \\
\hline
{4} & $N_4 = 0$ & {} & {} & {} & {} \\
\hline
\end{tabular}
\caption{\it Example on how the cube root of a thirteen digit number $N$ can be found in four iterations when the correct value for $a$ is selected. The algorithm terminates with a solution when $frac_i = f_i$. In this example $p = 2 + 9*5 + 9^2*3 + 9^3*1 = 1019$; and the cube root of $N$ is therefore: $7+18*1019 = 18349$.}
\end{center}
\end{figure}

\begin{figure}[h]\label{fig:in8}
\begin{center}
\begin{tabular}{|l|c|c|c|c|c|}
\hline
mod 18 & {} & {} & {[17]} & {[17]} & {[17]}\\ 
\hline
\hline
Twin Primes & $[r_i]$ & $[r_j]$ & Type A & Type B & Type C\\
\hline
\hline
{5, 7} & [5] & [7] & A & {} & {}\\
\hline
{11, 13} & [11] & [13] & {} & B & {}\\
\hline
{17, 19} & [17] & [1] & {} & {} & C\\
\hline
{29, 31} & [11] & [13] & {} & B & {}\\
\hline
41, 43 & [5] & [7] & A & {} & {}\\
\hline
{59, 61} & [5] & [7] & A & {} & {}\\
\hline
{71, 73} & [17] & [1] & {} & {} & C\\
\hline
{101, 103} & [11] & [13] & {} & B & {}\\
\hline
{107, 109} & [17] & [1] & {} & {} & C\\
\hline
{137, 139}  & [11] & [13] & {} & B & {}\\
\hline
{149, 151} & [5] & [7] & A & {} & {}\\
\hline
{179, 181} & [17] & [1] & {} & {} & C\\
\hline
\end{tabular}
\caption{\it Illustration on how twin primes can be divided into three different types based on their residues classes $(mod 18)$. Note how the product of all three types can be written as $17 + 18*X$.}
\end{center}
\end{figure}

\section{Twin Primes}
Since $U_{18}$ contains all odd primes, except 3, all twin primes are members of $U_{18}$. In fact, since the difference between two twin primes is $2$, twin primes can only occur between the residue classes [5], [7] (Type A); [11], [13] (Type B); and [17], [1] (Type C). The first 12 twin prime sets and their types are given in figure 8.

Because the product between twin primes always belongs to $[17]$ and squares always belong to $[1], [7], [13]$; the $a$, and therefore DR(N),  in $N = a + 18p$ becomes a convenient differentiator between the two cases.

\section{Acknowledgment}
The author would like to thank John Lindgren for providing very helpful suggestions, verifying essential elements, and most of all for encouraging the publication of this work.

\section{References}
\begin{enumerate}
\item Kenneth H. Rosen, Elementary Number Theory and its applications, Fourth Edition; May 2000.
\item G.H. Hardy and E.M. Wright, An introduction to the Theory of Numbers,Sixth Edition; September 2007.
\item US Patent Office Website; www.uspto.gov.
\item Victor Shoup, A Computational Introduction to Number Theory and Algebra, version 2.0; December 2007.
\end{enumerate}
\end{document}